\documentclass[preprint,12pt]{elsarticle}
\usepackage{amssymb}
\usepackage{amsmath}
\usepackage{amsfonts}
\begin{document}

\begin{frontmatter}

\title{Generalized even and odd totally positive matrices}
\author{O.Y. Kushel}
\ead{kushel@math.tu-berlin.de}
\address{Institut f\"{u}r Mathematik, MA 4-5, Technische Universit\"{a}t Berlin, D-10623 Berlin, Germany}
\author{P. Sharma}
\ead{sharmapatanjali@rediffmail.com}
\address{Dept. of Mathematics, University of Rajasthan, Jaipur 302 004, INDIA}
\begin{abstract}
A generalization of the definition of an oscillatory matrix based on the theory of cones is given in this paper. The positivity and simplicity of all the eigenvalues of a generalized oscillatory matrix are proved. The classes of generalized even and odd oscillatory matrices are introduced. Spectral properties of the obtained matrices are studied. Criteria of generalized even and odd oscillation are given. Examples of generalized even and odd oscillatory matrices are presented.
\end{abstract}


\end{frontmatter}

\newtheorem{thm}{Theorem}
\newtheorem{lem}[thm]{Lemma}
\newdefinition{rmk}{Remark}
\newproof{pf}{Proof}
\newproof{pot}{Proof of Theorem \ref{thm2}}
%%
%% Start line numbering here if you want
%%
% \linenumbers

%% main text
\section{Introduction}
The systematic theory of linear integral equations with totally positive and oscillatory kernels was developed mainly by F.R. Gantmacher and M.G. Krein in monograph \cite{GANT} and by S. Karlin in \cite{KARL}.
The following definition was given in \cite{GANT}:
the kernel $k(t,s) \in C[0,1]^2$ is called {\it oscillatory}, if it satisfies both the conditions:

(a) for any $0< t_0 < t_1 < \ldots < t_n < 1$ and $0< s_0 < s_1 <
\ldots < s_n < 1$, $ \ \ n = 0, \ 1, \ \ldots$ the inequality
$$k\begin{pmatrix}
  t_0 & t_1 & \ldots & t_n \\
  s_0 & s_1 & \ldots & s_n
\end{pmatrix} \geq 0$$
is true;

 (b) for any $0< t_0 < t_1 < \ldots < t_n <
1$, $ \ \ n = 0, \ 1, \ \ldots$ the inequality
 $$k\begin{pmatrix} t_0 & t_1 & \ldots & t_n \\
  t_0 & t_1 & \ldots & t_n
\end{pmatrix} > 0$$ is true.

 Such inequalities for the kernel allow to obtain remarkable spectral properties of the corresponding linear integral operator, that its eigenvalues are always simple, distinct and positive and its eigenunctions ocsillate in a definite manner.

The theory of total positivity is widely used for obtaining oscillation theorems for various boundary value problems, cause sometimes it is possible to show, that the Green's function of a considered differential equation is an oscillatory kernel.
However, futher studying of ordinary differential equations required certain generalization of the existing terminology. So, in papers \cite{KAL1} and \cite{KAL2} by P. Kalafati, devoted to studying boundary-value problems with quasi-differential operators, the more general definitions of even and odd oscillatory kernels were introduced.
The kernel $k(t,s)$ is called {\it even (odd) oscillatory}, if conditions (a) and (b) are true for all even (respectively odd) values of $n$.

Such a generalization gave an opportunity to get some information about the spectrum of a boundary-value problem in the case, when it is impossible to establish total positivity of its Green's function.

In turn, the aim of this paper is to generalize the conception of even and odd total oscillation, introduced by P. Kalafati, using the methods of the theory of cones. The authors believe that the developed apparatus would be suitable for applying to a wider class of boundary-value problem.
However, in this paper we shall be restricted to the simpliest finite-dimensional case.

\section{Exterior powers of the space ${\mathbb R}^n$}
First let us recall some definitions and statements of the exterior algebra.

Let $(i_1,  \ldots,  i_j)$ be a permutation of the indices $1,
\ldots,  j$. Define
$$\chi(i_1, \ldots, i_j) = \left\{\begin{array}{cc} 1,& \mbox{if
the permutation $(i_1, \ldots,  i_j)$ is even};
\\[10pt] - 1, & \mbox{if the permutation is odd.}\end{array}\right.$$

 Let $x_1, \ \ldots, \ x_j$ $ \ (2 \leq j \leq n)$ be arbitrary
vectors from $n$-dimensional Euclidean space ${\mathbb R}^n$. Then the multilinear functional $x_1
\wedge \ldots \wedge x_j: \times^j ({\mathbb R}^n)^ \prime
\rightarrow {\mathbb R}$, which acts according to the rule
$$(x_1 \wedge \ldots \wedge x_j)(f_1, \ldots, f_j)=
\sum_{(i_{1},  \ldots,  i_j)} \chi(i_1, \ldots, i_j) \ \langle
x_{i_1}, f_1 \rangle \ldots \langle x_{i_j}, f_j \rangle$$ is
called an {\it exterior product} of the vectors $x_1, \ \ldots, \
x_j$. Here the sum is taken with respect to all the permutations
$(i_{1},  \ldots,  i_j)$ of the set of indices $(1,  \ldots,
j)$, and linear functionals  $f_1, \ \ldots, \ f_n \in ({\mathbb
R}^n)^ \prime$ are considered as vectors from ${\mathbb R}^n$.

 The $j$th exterior power $\wedge^j {\mathbb R}^n$ is defined as a
linear span of all the exterior products of the form $x_1 \wedge
\ldots \wedge x_j$, where $x_1, \ \ldots, \ x_j \in {\mathbb
R}^n$. The space $\wedge^j {\mathbb R}^n$ is also finite
dimensional and is isomorphic to the space ${\mathbb R}^{C _n^j},$
where $C_n^j = \frac{n!}{j!(n-j)!}$. Examine an arbitrary basis
$e_1, \ \ldots, \ e_n$ in ${\mathbb R}^n$. The set of all exterior
products of the form $e_{i_1} \wedge \ldots \wedge e_{i_j}$, where
$1 \leq i_1 < \ldots < i_j \leq n$ forms a canonical basis in the
space $\wedge^j {\mathbb R}^n$ (see \cite{GLALU}).

Let $A$ be a linear operator in the space ${\mathbb R}^n$. Then the linear operator $\wedge^j A$, which acts in the space $\wedge^j
{\mathbb R}^n$ according to the rule:
$$ (\wedge^j A)(x_1 \wedge \ldots \wedge x_j) = Ax_1 \wedge \ldots \wedge Ax_j,$$
is called the {\it $j$th exterior power} of the operator $A$. It is easy to see, that $\wedge^1 A = A$ and $\wedge ^n A$ is one-dimensional and coinside with $\det A$.

If the operator $A$ is given by the matrix ${\mathbf A} =
\{a_{ij}\}_{i,j = 1}^n$ in the basis $e_1, \ \ldots, \ e_n$, then
the matrix of its $j$th exterior power $\wedge^j A$ in the basis
$\{e_{i_1} \wedge \ldots \wedge e_{i_j}\}$, where $1 \leq i_1<
\ldots < i_j \leq n$, coincides with the $j$th compound matrix $
{\mathbf A}^{(j)}$ of the initial matrix ${\mathbf A}$. Here the $j$th compound matrix $
{\mathbf A}^{(j)}$ consists of all the minors of the $j$th order
$A\begin{pmatrix}
  i_1 &  \ldots & i_j \\
  k_1 & \ldots & k_j \end{pmatrix}$, where $1 \leq i_1<
\ldots < i_j \leq n, \ 1 \leq k_1< \ldots < k_j \leq n$, of the
initial $n \times n$ matrix ${\mathbf A}$, numerated in the lexicographic order (see, for example,
\cite{PINK}).

Recall the statement, concerning the eigenvalues of the exterior
power of an operator, which we shall use later.

\begin{thm}[Kronecker] Let $\{\lambda_{i}\}_{i = 1}^n$ be the set of
all eigenvalues of a linear operator $A: {\mathbb R}^n \rightarrow
{\mathbb R}^n$, repeated according to multiplicity. Then all the
possible products of the form $\{\lambda_{i_1} \ldots \lambda_{i_j} \}$, where $1
\leq i_1 < \ldots < i_j \leq n$, forms the set of all the possible eigenvalues
of the $j$th exterior power $\wedge^j A$ of the operator $A$, repeated
according to multiplicity.
\end{thm}

The Kronecker theorem is formulated in terms of compound matrices
and proved without using exterior products in \cite{GANT} (see
\cite{GANT}, p. 80, theorem 23).

\section{Generalized even and odd oscillatory operators}
The classical methods of studying totally positive and oscillatory
matrices are based on the spectral analysis of nonnegative and primitive matrices. However, it was shown, (see \cite{BIR, KRU}), that the remarkable properties of the spectrum of nonnegative matrices also hold for operators, which leave
invariant a certain cone in ${\mathbb R}^n$. Let us recall some
definitions and statements of the theory of cone-preserving maps
(see, for example, \cite{BERPL}, \cite{TAM}).

A set $K \subset {\mathbb R}^n$ is called {\it a proper cone}, if
it is a convex cone (i.e. for any $x, y \in K, \ \alpha \geq 0$ we
have $x+y, \ \alpha x \in K$), is pointed (i.e. $K \cap (-K) =
\{0\}$), closed and solid (i.e. ${\rm int}(K) \neq \emptyset$).

Let $K \subset {\mathbb R}^n$ be a proper cone. A linear operator
$A:{\mathbb{R}}^{n} \rightarrow {\mathbb{R}}^{n}$ is called {\it
$K$-positive} or {\it positive with respect to the cone $K$}, if
$A(K \setminus \{0\}) \subseteq {\rm int}(K)$. In the case of $K =
{\mathbb{R}}^{n}_+$, $K$-positive operators are called simply {\it
positive}.

Let us formulate the following generalization of the Perron--Frobenius
theorem (see \cite{BERPL}, p. 13, theorem 3.26).

\begin{thm}[Generalized Perron] Let a linear operator $A : {\mathbb R}^n
\rightarrow {\mathbb R}^n$ be positive with respect to a proper
cone $K \subset {\mathbb R}^n$. Then:
\begin{enumerate}
\item[\rm 1.] The spectral radius $\rho(A)$ is a simple positive
eigenvalue of the operator $A$, different in modulus from the
other eigenvalues.
\item[\rm 2.] The eigenvector $x_1$, corresponding to
the eigenvalue $\lambda_1 = \rho(A)$, belongs to ${\rm int}(K)$.
\item[\rm 3.] The functional $x_1^*$, corresponding to
the eigenvalue $\lambda_1 = \rho(A)$, satisfy the inequality
$x_1^*(x) > 0$ for every nonzero $x \in K$.
\end{enumerate}
\end{thm}

Later we shall also deal with the following generalization of the definition of
positivity.

A linear operator $A: {\mathbb R}^n \rightarrow {\mathbb
R}^n$ is called {\it $K$--primitive} or {\it primitive with respect to the cone $K$}, if $AK \subseteq K$ and the only nonempty
subset of $\partial(K)$, which is left invariant by $A$, is $\{0\}$.
This definition was given by Barker (see \cite{BARK}, see also
\cite{TAM}).

The following statement is true for $K$-primitive operators:
{\it a linear operator $A$ is primitive with respect to some proper cone $K$ if and only if there exists such a proper cone $\widetilde{K}$, that $A$ is positive with respect to $\widetilde{K}$} (see \cite{BERPL}, p. 18, corollary 4.13).
 So if the operator $A$ is $K$-primitive with respect to some
proper cone $K$, then the statement of Theorem 2 is true for
the operator $A$ as well. In general, the study of $K$-primitive operators amounts to the study of $K$-positive operators, but in some cases it is much easier to prove the primitivity with respect to a certain cone in ${\mathbb R}^n$, then the positivity with respect to some other cone.

To introduce the generalization of the class of oscillatory
operators, we shall use the method of the crossway from the study
of the initial operator in the space ${\mathbb R}^n$ to the study of its exterior powers $ \wedge^j A$, $(j = 2, \ \ldots, \ n)$, acting, respectively, in the exterior powers of the initial space $\wedge^j {\mathbb R}^n = {\mathbb R}^{C_n^j}$. Let us give the following definitions.

A linear operator $A$ is called {\it generalized oscillatory (GO)} if
it is primitive with respect to a proper cone $K_1 \subset
{\mathbb R}^n$, and for every $j$ $(j = 2, \ \ldots, \ n)$ its $j$-th
exterior power $ \wedge^j A$ is also primitive with respect to a
proper cone $K_j \subset {\mathbb R}^{C_n^j}$.

In the case, when $K_1 = {\mathbb R}^n_+$ and $K_j = {\mathbb
R}^{C_n^j}_+$ for every $j$ $(j = 2, \ \ldots, \ n)$, the definition, given above, coincides with the classical definition of an oscillatory matrix,
given by F.R. Gantmacher and M.G. Krein in \cite{GANT}.

Let us notice, that every matrix, similar to some oscillatory one, is generalized oscillatory. But the inverse statement is wrong: not every generalized oscillatory matrix is similar to some oscillatory matrix.

In turn, the definitions of even (odd) oscillation will be generalized in the following way.

A linear operator $A$ is called {\it generalized even
oscillatory (GEO) (or generalized odd oscillatory (GOO))} if for every even (respectively odd) $j$ $(j = 1, \ \ldots, \ n)$ its $j$-th exterior power $ \wedge^j A$ is primitive
with respect to a proper cone $K_j \subset {\mathbb R}^{C_n^j}$.

\section{Spectral properties of generalized even and odd oscillatory operators}
Let us formulate and prove the following theorems about spectral properties of GO operators.

\begin{thm} Let a linear operator $A: {\mathbb R}^n
\rightarrow {\mathbb R}^n$ be generalized oscillatory. Then all
the eigenvalues of the operator $A$ are simple, positive and
different in modulus from each other:
$$\rho(A) = \lambda_1 > \lambda_2 > \ldots > \lambda_n > 0.$$ \end{thm}

\begin{pf}
Enumerate the eigenvalues of the operator $A$
in order of decrease of their modules (taking into account their
multiplicities):
$$|\lambda_{1}| \geq | \lambda_{2}| \geq |\lambda_{3}| \geq \ldots
\geq |\lambda_{n}|.$$
 Applying the generalized Perron theorem to the operator $A$,
we get: $\lambda_{1} = \rho(A)>0$ is a simple positive eigenvalue
of $A$, different in modulus from the rest of eigenvalues. Examine the second exterior power $\wedge^2 A$, which is also positive with respect to some proper cone $K_2 \subset {\mathbb R}^{C_n^2}$.
Applying generalized Perron theorem to $\wedge^2 A$, we get: $\rho(\wedge^2 A) > 0$ is
also a simple positive eigenvalue of $\wedge^2 A$, different in modulus from the rest of eigenvalues.

 As it follows from the statement of Kroneker theorem, $\wedge^2 A$ has no other eigenvalues, except all the possible products of the form $\lambda_{i_1}\lambda_{i_2}$, where $1 \leq i_1 < i_2 \leq n$. Therefore $\rho(\wedge^2 A)>0$ can be represented in the
form of the product $\lambda_{i_1}\lambda_{i_2}$ with some values of
the indices $i_1,i_2$, \ $i_1 < i_2$. It follows from the facts that the
eigenvalues are numbered in a decreasing order, and there is only
one eigenvalue on the spectral circle $|\lambda| = \rho(A)$, that
 $\rho(\wedge^2 A) =  \lambda_{1}\lambda_{2} = \rho(A)\lambda_{2}$.
  Therefore $\lambda_{2} = \frac{\rho(\wedge^2 A)}{\rho(A)}>0$.

Repeating the above reasoning for $\wedge^j A$, $j = 3, \ \ldots, \ n$, we receive the relations:
$$\lambda_{j} = \frac{\rho(\wedge^j A)}{\rho(\wedge^{j-1} A)}>0,$$
where $j = 3, \ \ldots, \ n$. The simplisity and distinction of the eigenvalues $\lambda_j$ for every $j$ also follows from the above relations, and the simplisity and distinction of $\rho(\wedge^j A)$.
\end{pf}

A weaker statements are true for GEO and GOO operators.

\begin{thm} Let a linear operator $A: {\mathbb R}^n
\rightarrow {\mathbb R}^n$ be even generalized oscillatory. Then
the algebraic multiplicity $m(\lambda)$ of any eigenvalue
$\lambda$ of the operator $A$ is not greater than $2$. The
following inequalities for the modules of the eigenvalues are
true:
$$\rho(A) = |\lambda_1| \geq |\lambda_2| > |\lambda_3| \geq |\lambda_4| > \ldots. $$
(The eigenvalues of $A$ are repeated according to multiplicity in
the above numeration.) Moreover, for every pair
$\lambda_i\lambda_{i+1} \ (i = 1, \ 3, \ 5, \ \ldots)$ the
following equality is true: $\arg(\lambda_{i+1}) = -
arg(\lambda_{i})$. If $n$ is odd, then $\lambda_n$ is real.\end{thm}

\begin{pf}
Enumerate the eigenvalues of the operator $A$
in order of decrease of their modules (taking into account their
multiplicities):
$$|\lambda_{1}| \geq | \lambda_{2}| \geq |\lambda_{3}| \geq \ldots
\geq |\lambda_{n}|.$$

Applying the generalized Perron theorem to the operator $\wedge^2 A$ we get:
$\rho(\wedge^2 A) > 0$ is a simple positive eigenvalue of $\wedge^2 A$, different in modulus from the rest of eigenvalues.

Let us prove, that there are no more than $2$ eigenvalues on the spectral circle $|\lambda| = \rho(A)$. Assume the contrary: there are $m > 2$ eigenvalues on the largest spectral circle. Let $M = \{1, \ \ldots, \ m\}$ be the set of their numbers. As it follows from the Kroneker theorem, $\rho(\wedge^2 A)>0$ can be represented
in the form of the product $\lambda_{i_1}\lambda_{i_2}$ with some
values of the indices $i_1,i_2$, \ $1 \leq i_1 < i_2 \leq n$. Since the eigenvalues are numerated in the decreasing order, then the equality $ \rho(\wedge^2 A) = |\lambda_{i_1}\lambda_{i_2}|$ is true for every pair $i_1,i_2$, where $i_1, \ i_2 \in M$ and $i_1 < i_2$. We have received, that $\rho(\wedge^2 A)$ is either a multiple eigenvalue of $\wedge^2 A$ or there are other eigenvalues, equal in modulus to $\rho(\wedge^2 A)$. This contradicts the generalized Perron theorem.
 So we have got the inequality
$$\rho(A) = |\lambda_1| \geq |\lambda_2| > |\lambda_3| \geq |\lambda_4| > \ldots. $$
As it follows, we can represent $\rho(\wedge^2 A) = \lambda_1\lambda_2$.
 In turn, it follows from the positivity of $\rho(\wedge^2 A)$, that $\lambda_{1}$ and $\lambda_2$ are either
a pair of complex adjoint eigenvalues, or both are real and of the same sign.

Repeating the above reasoning for the operator $\wedge^j A$, $j = 4, \ 6, \ \ldots, \ 2[\frac{n}{2}]$, we receive the relations:
$$\rho(\wedge^j A) = \prod_{i =1}^j \lambda_i = \rho(\wedge^{j-2} A)\lambda_{j-1}\lambda_j.$$
So it follows the simplisity and distinction of $\rho(\wedge^j A)$, that there is at most two eigenvalues on every spectral circle of the operator $A$.
The reality of the last eigenvalue $\lambda_{n}$ easily follows from the equality
$$\lambda_{n} = \frac{\det A}{\prod_{i =1}^{n-1} \lambda_i} = \frac{\det A}{\rho(\wedge^{n-1} A)},$$  in the case when $n-1$ is even.
\end{pf}

\begin{thm} Let a linear operator $A: {\mathbb R}^n
\rightarrow {\mathbb R}^n$ be odd generalized oscillatory. Then
the algebraic multiplicity $m(\lambda)$ of any eigenvalue
$\lambda$ of the operator $A$ is not greater than $2$. The
following inequalities for the modules of the eigenvalues are
true:
$$\rho(A) = |\lambda_1| > |\lambda_2| \geq |\lambda_3| > |\lambda_4| \geq \ldots. $$
(The eigenvalues of $A$ are repeated according to multiplicity in
the above numeration.) Moreover, $\lambda_1 = \rho(A)$ is a simple
positive eigenvalue of $A$. If $n$ is even, then $\lambda_n$ is
real. For every pair $\lambda_i\lambda_{i+1} \ (i = 2, \ 4, \ 6, \
\ldots)$ the following equality is true: $\arg(\lambda_{i+1}) = -
arg(\lambda_{i})$.\end{thm}
\begin{pf}
The proof is same as the proof of Theorem 4.
\end{pf}

\section{Some matrix criteria for generalized even and odd oscillation}
Let us reformulate the above theorems in terms of compound
matrices. Then the conditions of this theorems will become easily
verified. First let us remind some well-known definitions and
statements.

 A matrix ${\mathbf A}$ is called
{\it non-negative} ({\it positive}), if all its elements $a_{ij}$
are nonnegative (positive). A nonnegative matrix ${\mathbf A}$ is
called {\it primitive} if there exists such a natural number $m$,
that the matrix ${\mathbf A}^m$ is positive. If the matrix
${\mathbf A}$ of a linear operator $A: {\mathbb R}^n \rightarrow
{\mathbb R}^n$ is primitive, then $A$ is $K$-primitive with
respect to the cone $K_+$ of all nonnegative vectors from the
space ${\mathbb R}^n$. The statement, that if the matrix
$\mathbf{A}$ is similar to a primitive matrix, then the
corresponding operator $A$ is $K$-primitive with respect to some
polyhedral cone $K$ in ${\mathbb R}^n$, easily follows from the
above reasoning. In some special cases we can see, if the matrix
${\mathbf A}$ is similar to a primitive matrix, just looking at
its structure. Let us recall the following two definitions (see, for example, \cite{KU1}).

A matrix $\mathbf A$ is called {\it strictly ${\mathcal
J}$--sign-symmetric}, if ${\mathbf A}$ does not contain zero
elements and there exists such a subset ${\mathcal J} \subseteq
\{1, \ \ldots, \ n\}$, that the inequality $a_{ij} < 0$ is true if
and only if one of the numbers $i$, $j$ belongs to the set
${\mathcal J}$, and the other belongs to the set $\{1, \ \ldots, \
n\}\setminus {\mathcal J}$.

A matrix ${\mathbf A}$ of a linear operator $A: {\mathbb R}^n
\rightarrow {\mathbb R}^n$ is called {\it ${\mathcal
J}$--sign-symmetric}, if there exists such a subset ${\mathcal J}
\subseteq \{1, \ \ldots, \ n\}$, that both the conditions (a) and
(b) are true:
\begin{enumerate}
\item[\rm (a)] the inequality $a_{ij} \leq 0$ follows from the inclusions $i
\in {\mathcal J}$, $j \in \{1, \ \ldots, \ n\}\setminus {\mathcal
J}$ and from the inclusions $j \in {\mathcal J}$, $i \in \{1, \
\ldots, \ n\}\setminus {\mathcal J}$ for any two numbers $i,j$;

\item[\rm (b)] one of the inclusions $i \in {\mathcal J}$, $j \in \{1, \
\ldots, \ n\}\setminus {\mathcal J}$ or $j \in {\mathcal J}$, $i
\in \{1, \ \ldots, \ n\}\setminus {\mathcal J}$ follows from the
strict inequality $a_{ij} < 0$.
\end{enumerate}

If the matrix ${\mathbf A}$ is ${\mathcal J}$--sign-symmetric
(strictly ${\mathcal J}$--sign-symmetric), then it is similar to
some nonnegative (respectively positive) matrix. Moreover, the
matrix of the similarity transformation is diagonal, and its
diagonal elements are equal to $\pm 1$ (see \cite{KU1}).

It's easy to see, that if the matrix ${\mathbf A}$ is ${\mathcal
J}$--sign-symmetric, and the matrix ${\mathbf A}^m$ is strictly
${\mathcal J}$--sign-symmetric for some natural number $m$, then
the matrix ${\mathbf A}$ is similar to some primitive matrix with
the diagonal matrix of the similarity transformation. Let us call
such matrices {\it ${\mathcal J}$--sign-symmetric primitive}. In
this case the linear operator $A: {\mathbb R}^n \rightarrow
{\mathbb R}^n$, defined by the matrix ${\mathbf A}$, is
$K$-primitive with respect to some cone spanned by the vectors
$e'_1, \ \ldots, \  e'_n$, where each vector $e'_i$ is equal
either to $e_i$ or to $- e_i$ $ \ (i = 1, \ \ldots, \ n)$.

Using the above reasoning, we receive the following theorems.

\begin{thm} Let the matrix ${\mathbf A}$ of a linear
operator $A: {\mathbb R}^n \rightarrow {\mathbb R}^n$ be
${\mathcal J}$--sign-symmetric primitive, and let the $j$th
compound matrix ${\mathbf A}^{(j)}$ be also ${\mathcal
J}$--sign-symmetric primitive for every $j$ $(1 < j \leq n)$. Then
all the eigenvalues of the operator $A$ are simple, positive and
different in modulus from each other:
$$\rho(A) = \lambda_1 > \lambda_2 > \ldots > \lambda_n > 0.$$ \end{thm}

\begin{thm} Let the $j$th compound matrix ${\mathbf
A}^{(j)}$ of the matrix ${\mathbf A}$ of a linear operator $A:
{\mathbb R}^n \rightarrow {\mathbb R}^n$ be ${\mathcal
J}$--sign-symmetric primitive for every even $j$ $(1 \leq j \leq
n)$. Then the algebraic multiplicity $m(\lambda)$ of any
eigenvalue $\lambda$ of the operator $A$ is not greater than $2$.
The following inequalities for the modules of the eigenvalues are
true:
$$\rho(A) = |\lambda_1| \geq |\lambda_2| > |\lambda_3| \geq |\lambda_4| > \ldots. $$
(The eigenvalues of $A$ are repeated according to multiplicity in
the above numeration.) Moreover, for every pair
$\lambda_i\lambda_{i+1} \ (i = 1, \ 3, \ 5, \ \ldots)$ the
following equality is true: $\arg(\lambda_{i+1}) = -
arg(\lambda_{i})$. If $n$ is odd, then $\lambda_n$ is real. \end{thm}

\begin{thm} Let the $j$th compound matrix ${\mathbf
A}^{(j)}$ of the matrix ${\mathbf A}$ of a linear operator $A:
{\mathbb R}^n \rightarrow {\mathbb R}^n$ be ${\mathcal
J}$--sign-symmetric primitive for every odd $j$ $(1 \leq j \leq
n)$. Then the algebraic multiplicity $m(\lambda)$ of any
eigenvalue $\lambda$ of the operator $A$ is not greater than $2$.
The following inequalities for the modules of the eigenvalues are
true:
$$\rho(A) = |\lambda_1| > |\lambda_2| \geq |\lambda_3| > |\lambda_4| \geq \ldots. $$
(The eigenvalues of $A$ are repeated according to multiplicity in
the above numeration.) Moreover, $\lambda_1 = \rho(A)$ is a simple
positive eigenvalue of $A$. If $n$ is even, then $\lambda_n$ is
real. For every pair $\lambda_i\lambda_{i+1} \ (i = 2, \ 4, \ 6, \
\ldots)$ the following equality is true: $\arg(\lambda_{i+1}) = -
arg(\lambda_{i})$.\end{thm}

\section{Examples}

Let us give some examples, illustrating the above theorems.

{\bf Example 1.}
Let the operator $A:{\mathbb{R}}^{3} \rightarrow
{\mathbb{R}}^{3}$ be defined by the matrix $$\mathbf{A} =
\begin{pmatrix}
  4 & -6.8 & 4.4 \\
 -1.2 & 6.3 & -1.1 \\
  1.8 & -2.6 & 3.4 \\
\end{pmatrix}.$$

 This matrix is strictly ${\mathcal J}$--sign-symmetric.
In this case the set ${\mathcal J}$ in the definition of
${\mathcal J}$--sign-symmetricity of the matrix $\mathbf{A}$
consists of two indices $1$ and $3$.

In this case the second compound matrix is the following:
$${\mathbf A}^{(2)} = \begin{pmatrix}
 17.04 & 0.88 & -20.24 \\
 1.84 & 5.68 & -11.68 \\
  -8.22 & -2.1 & 18.56 \\
\end{pmatrix}.$$

The matrix ${\mathbf A}^{(2)}$ is also strictly
${\mathcal J}$--sign-symmetric. The set
the set ${\mathcal J}$ in the definition of
${\mathcal J}$--sign-symmetricity of the matrix ${\mathbf A}^{(2)}$
consists of two indices $1$ and $2$.

The third compound matrix ${\mathbf A}^{(3)}$ consists of only one positive element, which is equal to $23.792$, so it can also be considered as strictly
${\mathcal J}$--sign-symmetric.

The operator $A$ satisfies the
conditions of theorem 6. It is easy to see, that all the three eigenvalues of the operator $A$, which are $\lambda_1 = 9.69542$, $\lambda_2 = 3.24937$, $ \lambda_3 = 0.755205$, are simple, positive and different in modulus from each other.

{\bf Example 2.}
Let the operator $A:{\mathbb{R}}^{4} \rightarrow
{\mathbb{R}}^{4}$ be defined by the matrix $$\mathbf{A} =
\begin{pmatrix}
 7 & 5.2 & 7.8 & 18.6 \\
 -6.9 & 4.4 & 5.3 & 37.5 \\
  2.1 & 4 & 5.6 & 20.8 \\
-9 & 1.8 & -2.4 & 17.4 \\
\end{pmatrix}.$$

In this case the second compound matrix is the following:
$${\mathbf A}^{(2)} = \begin{pmatrix}
 66.68 & 90.92 & 390.84 & -6.76 & 113.16 & 193.92 \\
 17.08 & 22.82 & 106.54 & -2.08 & 33.76 & 58.08 \\
 59.4 & 53.4 & 289.2 & -26.52 & 57 & 180.36 \\
-36.84 & -49.77 & -222.27 & 3.44 & -58.48 & -99.76 \\
27.18 & 64.26 & 217.44 & -20.1 & 9.06 & 182.22 \\
39.78 & 45.36 & 223.74 & -19.68 & 32.16 & 147.36 \\
\end{pmatrix}.$$

The matrix ${\mathbf A}^{(2)}$ is strictly
${\mathcal J}$--sign-symmetric. The set
the set ${\mathcal J}$ for the matrix ${\mathbf A}^{(2)}$
consists of only one index $4$.

The fours compound matrix ${\mathbf A}^{(4)}$ consists of only one positive element, which is equal to $278.964$, so it can also be considered as strictly
${\mathcal J}$--sign-symmetric.

The operator $A$ satisfies the
conditions of theorem 7 for GEO operators. It is easy to see, that the operator $A$
has two pairs of complex adjoint eigenvalues: $\lambda_1 = 17.813 + 16.2621 i$, $\lambda_2 = 17.813-16.2621 i$ on the first spectral circle and $\lambda_3 = -0.613045+0.322013 i$, $\lambda_4 = -0.613045-0.322013 i$ on the second spectral circle.

{\bf Example 3.}

Let the operator $A:{\mathbb{R}}^{4} \rightarrow
{\mathbb{R}}^{4}$ be defined by the matrix $$\mathbf{A} =
\begin{pmatrix}
1 & 8 & 3 & 0.4 \\
5.7 & 7.4 & 8.7 & 9.5 \\
1.5 & 9.7 & 2.5 & 6 \\
4 & 8.6 & 9.9 & 2.2 \\
\end{pmatrix},$$
which is obviously positive.

Its third compound matrix is the following:
$${\mathbf A}^{(3)} = \begin{pmatrix}
 57.08 & -189.674 & -30.92 & 344.494 \\
 -116.34 & 146.028 & 10.122 & -403.644 \\
 -41.97 & 124.98 & 10.14 & -310.608 \\
 163.601 & -265.352 & -81.065 & 572.437 \\
\end{pmatrix}.$$

The matrix ${\mathbf A}^{(3)}$ is strictly
${\mathcal J}$--sign-symmetric. The set
the set ${\mathcal J}$ for the matrix ${\mathbf A}^{(3)}$
is $\{2, \ 3\}$.

The operator $A$ satisfies the
conditions of theorem 8 for GOO operators. Then the operator $A$ has
the first positive simple eigenvalue $\lambda_1 = \rho(A) =
23.8704$ which is different in modulus from the other eigenvalues, a pair of complex adjoint eigenvalues on the second spectral circle $\lambda_{2} =
-5.58952 + 2.36837 i$,  $\lambda_3 = -5.58952 - 2.36837 i$, and the least in modulus real eigenvalue $\lambda_4 = 0.408632$.
%% References
%%
%% Following citation commands can be used in the body text:
%% Usage of \cite is as follows:
%%   \cite{key}         ==>>  [#]
%%   \cite[chap. 2]{key} ==>> [#, chap. 2]
%%

%% References with bibTeX database:

\bibliographystyle{elsarticle-num}
\bibliography{<your-bib-database>}

%% Authors are advised to submit their bibtex database files. They are
%% requested to list a bibtex style file in the manuscript if they do
%% not want to use elsarticle-num.bst.

%% References without bibTeX database:

% \begin{thebibliography}{00}

%% \bibitem must have the following form:
%%   \bibitem{key}...
%%

% \bibitem{}

% \end{thebibliography}

\end{document}